\documentclass[
leqno,          
11pt,           
a4paper         
]{amsart}

\usepackage[colorlinks, citecolor=blue,pagebackref,hypertexnames=false
]{hyperref}
\usepackage{courier}                                                    
\usepackage{amssymb, amsmath, amsfonts, amsthm}                         
\usepackage{mathtools, cite, enumerate, rotating, framed, lipsum,enumitem,esint}       
\usepackage{fancybox, bbm}                                              
\usepackage[dvipsnames]{xcolor}                                         
\usepackage[polish, english]{babel}                                     
\usepackage[T1]{fontenc}                                                
\usepackage[cp1250]{inputenc}                                           
\usepackage[normalem]{ulem}
\usepackage{yfonts}
\usepackage{mathrsfs}

\DeclareMathAlphabet{\mathpzc}{OT1}{pzc}{m}{it}
 \numberwithin{thmcount}{section}   
 \numberwithin{equation}{section}                        

\newcommand{\thmcount}{thmcount}                 
\newcounter{specialcounter}

\newtheorem{Thm}[\thmcount]{Theorem}
\newtheorem{Sthm}[specialcounter]{Theorem}
\newtheorem{Cor}[\thmcount]{Corollary}
\newtheorem{Lem}[\thmcount]{Lemma}
\newtheorem{Prop}[\thmcount]{Proposition}
\newtheorem{Rem}[\thmcount]{Remark}
\newtheorem{Defn}[\thmcount]{Definition}
\newtheorem{Ex}[\thmcount]{Example}
\newtheorem{Asu}[\thmcount]{Assumption}
\newtheorem{Sol}[\thmcount]{Solution}
\newtheorem*{Thmx}{Theorem}
\newtheorem*{Corx}{Corollary}
\newtheorem*{Lemx}{Lemma}
\newtheorem*{Propx}{Proposition}

\newtheorem*{Defnx}{Definition}
\newtheorem*{Exx}{Example}
\newtheorem*{Asux}{Assumption}
\newtheorem*{Solx}{Solution}

\newcommand \eq[1]{\begin{equation} #1 \end{equation}}

\renewcommand \sp[1]{\begin{equation} \begin{split} #1 \end{split} \end{equation}}
\newcommand \spx[1]{\begin{equation*} \begin{split} #1 \end{split} \end{equation*}}
\newcommand \en[1]{\begin{enumerate}  #1 \end{enumerate}}

\newcommand{\thm}[2]{\begin{Thm} \label{#1} #2 \end{Thm}}
\newcommand{\sthm}[2]{\begin{Sthm} \label{#1} #2 \end{Sthm}}

\newcommand{\lem}[2]{\begin{Lem} \label{#1} #2 \end{Lem}}

\newcommand{\cor}[2]{\begin{Cor} \label{#1} #2 \end{Cor}}

\newcommand{\prop}[2]{\begin{Prop} \label{#1} #2 \end{Prop}}

\newcommand{\rem}[2]{\begin{Rem} \label{#1} #2 \end{Rem}}

\newcommand{\defn}[2]{\begin{Defn} \label{#1} #2 \end{Defn}}

\newcommand{\pr}[1]{\begin{proof} #1 \end{proof}}

\newcounter{comcount}

      \newcommand{\blue}{\textcolor{blue}}

        \newcommand{\la}{\lambda}

\newcommand{\CC}{\mathbb{C}}

\newcommand{\HH}{\mathbb{H}}  \newcommand{\hh}{\mathcal{H}}

  \renewcommand{\ll}{\mathcal{L}}
  
\newcommand{\NN}{\mathbb{N}}

\newcommand{\RR}{\mathbb{R}}

\newcommand{\supp}{\mathrm{supp}}
\newcommand{\8}{\infty}

\newcommand{\Rn}{{\RR^n}}

\newcommand{\wt}[1]{\widetilde{#1}}

\newcommand{\st}{\sqrt{t}}

\newcommand{\abs}[1]{\left| #1 \right|}
\newcommand{\set}[1]{\left\{ #1 \right\}}
\newcommand{\norm}[1]{\left\| #1 \right\|}

\newcommand{\eee}[1]{\left( #1 \right)}

\newcommand{\BMO}{\mathrm{BMO}(\mu; h_1,...,h_m;r')}
\newcommand{\Har}{H^1_{at}(\mu;h_1,...,h_m;r)}
\newcommand{\wh}{\widehat}

\begin{document}

\title{Hardy spaces meet harmonic weights revisited}

 \author[M. Preisner  and  A. Sikora]{ Marcin Preisner  \ and  Adam Sikora}
\address{Marcin Preisner, Instytut Matematyczny, Uniwersytet Wroc\l awski, \ pl. Grunwaldzki 2, 50-384 Wroc\l aw, Poland }
\email{marcin.preisner@uwr.edu.pl}
\address{
	Adam Sikora, Department of Mathematics, Macquarie University, NSW 2109, Australia}
\email{adam.sikora@mq.edu.au}

\begin{abstract}
We investigate  Hardy spaces $H^1_L(X)$ corresponding to self-adjoint operators~$L$. Our main aim is to obtain a description of $H^1_L(X)$ in terms of atomic decompositions similar to such characterisation of the classical Hardy spaces $H^1(\mathbb{R}^d)$. 
Under suitable assumptions, such a description was obtained by Yan and the authors in [Trans. Amer. Math. Soc. 375 (2022), no. 9, 6417-6451], where 
the atoms associated with an $L$-harmonic function are considered. Here we continue this study and modify the previous definition of atoms.

The modified approach allows us to investigate settings, when the generating operator is related to a system of linearly independent harmonic functions. In this context, the cancellation condition for atoms is adjusted to fit this system. In an explicit example, we consider a~symmetric manifold with ends $\RR^d \# \RR^d$. For this manifold the space of bounded harmonic functions is two-dimensional. 
Any element from the Hardy space $H^1_L(X)$ has to be orthogonal to all of the harmonic functions in the system.

\end{abstract}


\maketitle

{\it In memoriam of Jacek Zienkiewicz.}

\section{Introduction}

\subsection{Background and assumptions}\label{ssec11}

Hardy and BMO spaces theory is a significant mainstream area of harmonic analysis. Its origin can be traced back to a notion of complex analysis
initiated by Hardy. For an elegant and comprehensive description of the classical theory of real Hardy and BMO spaces, we refer readers to \cite{Stein}. A recent new direction in this area is to extend the investigation to Hardy and BMO spaces corresponding to semigroups generated by some differential operators. It is a natural idea. Indeed, one of the equivalent forms of the classical definition can be stated in terms of the maximal function corresponding to the heat semigroup. That is, we define the maximal function 
$M$ by the formula 
\begin{equation*}
Mf(x)=\sup_{t>0} |\exp(t\Delta)f(x)|, \quad x\in \Rn,
\end{equation*}
where $\Delta$ is the standard Laplace operator and $\exp(t\Delta)$ is the corresponding heat semigroup. Then we define the $H^1(\Rn)$-norm by setting 
\begin{equation*}
\|f\|_{H^1(\RR^n)}= \|Mf(x)\|_{L^1(\RR^n)}.
\end{equation*}

This definition can be generalised for any semigroup generated by a non-negative self-adjoint operator $L$ defined on $L^2(X,\mu)$, where $(X,\mu)$ is a measure space. There are many papers devoted to the study of Hardy and BMO spaces corresponding to the semigroup of operators, see for example \cite{DZ_JFAA, DZ_Revista2, Song-Yan-2017, Duong_Yan_2005, Chang_Krantz_Stein, Auscher_Russ, Hofmann_Mayboroda_McIntosh} and references therein. A general approach in this direction is discussed in the seminal paper \cite{Hofmann_Memoirs}. 

Our main goal is to describe the Hardy spaces $H^1_L(X)$ associated with an operator $L$ in terms of atomic decompositions. For the classical Hardy space $H^1(\Rn)$ such theory was obtained in \cite{Coifman_Studia, Latter_Studia}.  In \cite{Preisner_Sikora_Yan} the authors described the Hardy space $H^1_L(X)$ by using atoms that are related to a non-negative $L$-harmonic function $h$. Here we continue this study and modify the notion of atoms from there. The new definition is more natural and its geometrical aspects are more evident. The approach which we develop in this paper  allows us to study the cases, when there exist two or more linearly independent harmonic functions related to operator $L$. In this situation,  atoms satisfy cancellation conditions with respect to all these functions.  Additionally, we study dual BMO-type spaces.


We shall work under the following assumptions. Let $(X,d,\mu)$ be a metric measure space of homogeneous type in the sense of Coifman and Weiss.  In particular, $\mu$ satisfies the doubling condition and there exists $C, D>0$ such that
\eq{\label{doubling}
\mu(B(x,R))\leq C\eee{\frac{R}{r}}^D \mu(B(x,r)), \quad x\in X, \ 0<r\leq R,
}
where $B(x,r) = \set{y\in X \, : \, d(x,y) <r}$.

 In what follows we study a positive self-adjoint operator $L$ acting  on $L^2(X,\mu)=L^2(\mu)$ and the semigroup $T_t = \exp(-tL)$ generated by $L$. We assume that there exists a non-negative integral kernel $T_t(x,y)$ related to $T_t$ that satisfies the Gaussian upper bounds, i.e.
there exist $C,c>0$ such that
\eq{\tag{UG}\label{UG}
0 \leq T_t(x,y) \leq C \mu(B(x,\sqrt{t}))^{-1/2} \exp\eee{-\frac{d(x,y)^2}{ct}}, \quad t>0, \ x,y\in X.
}
The maximal function corresponding to the semigroup $T_t = \exp(-tL)$
is given by the formula 
$
M_Lf(x)=\sup_{t>0} \abs{T_tf(x)}.
$
For $f\in L^2(\mu)$ we set
\begin{equation*}
\|f\|_{H^1_L(X)}= \|M_Lf\|_{L^1(\mu)}.
\end{equation*}
Define the Hardy space $H^1_L(X)$ as the completion of the set $\{f\in L^2(\mu) \, : \, \|f\|_{H^1_L(X)} <\8\}$ with respect to this norm, c.f. \cite{Hofmann_Memoirs, Preisner_Sikora_Yan}. 

In this paper we shall consider the following atomic Hardy spaces.
\defn{def_atoms}{
Let $r,r'>1$ be the conjugate exponents, i.e. $1/r+1/r' = 1$, $m\in \NN$, and $h_1,...,h_m \in L^{r'}_{loc}(\mu)$. We say that a~function $a: \, X \to \CC$ is a $(\mu;h_1,...,h_m;r)$-atom  if there exists a ball $B$ in X such that:
\eq{\label{new_atoms}
\supp \, a \subseteq B, \qquad \norm{a}_{L^r\eee{\mu}} \leq \mu(B)^{-1+1/r}, \qquad \int a(x)h_j(x)\, d\mu(x)=0
}
for $j=1,...,m$. The space $\Har$ is defined as follows: we say that a function $f$ is in $\Har$ if $f(x) =\sum_k \la_k a_k(x)$, where $\sum_k |\la_k|<\8$ and $a_k$ are $(\mu;h_1,...,h_m;r)$-atoms. We set
$$\norm{f}_{\Har} = \inf \sum_k|\la_k|,$$
where the infimum is taken over all representations as above.
}

Note that for an atom $a$ as in \eqref{new_atoms} we have $\norm{a}_{L^1(\mu)}\leq 1$ and thus the series $\sum_k \la_k a_k(x)$ converges in $L^1(\mu)$. Moreover, by a standard argument, $H^1_{at}(\mu;h_1,...,h_m;r)$ is a~Banach subspace of $L^1(\mu)$. We shall also study related $BMO$-type spaces that are defined as follows.

\defn{BMO_general}{
Let $m\in \NN$, $r'>1$, and $h_1,...,h_m \in L^{r'}_{loc}(\mu)$.  Define $V = \mathrm{span}( h_1,...,h_m)$. 
For $g \in L^{r'}_{loc}(\mu)$ we define
$$\norm{g}_{BMO(\mu;h_1,...,h_m;r')} =\sup_{B} \eee{ \inf_{u \in V} \eee{\mu(B)^{-1}\int_B \abs{g(x) - u(x)}^{r'} d\mu(x)}^{1/r'}},$$
where the supremum is taken over all balls in $X$. By definition, the elements of the space  $BMO(\mu;h_1,...,h_m;r')$ are classes $\{g+u \, : \, u\in V\}$ such that $\norm{g}_{BMO(\mu;h_1,...,h_m;r')}<\8$.
}

\subsection{Operators related to one harmonic function} \label{ssec12}

In this section, as in \cite{Preisner_Sikora_Yan}, we assume that there exists a non-negative harmonic function $h$ such that the heat semigroup $T_t$, after the Doob transform (or $h$-transform) corresponding to $h$, satisfies the upper and lower Gaussian estimates. More precisely, we assume that there exists a function $h\in L^\8_{loc}(X)$ \blue{ that is positive a.e., $h$} is $L$-harmonic in the sense that for all $t>0$
\eq{\label{L-harm}
T_t h (x) =h(x), \quad \text{a.e. } x\in X,
}
and that there exist $c_1, c_2, C>0$ such that for a.e. $x,y\in X$ and $t>0$ we have
 \begin{equation}
 \label{ULG'}\tag{ULG\textsubscript{h}}
 \frac{C^{-1}}{\mu_{h^2}(B(x,\st))} \exp\eee{-\frac{d(x,y)^2}{c_1t}}\leq \frac{T_t(x,y)}{h(x)h(y)}
 \leq \frac{C}{\mu_{h^2}(B(x,\st))} \exp\eee{-\frac{d(x,y)^2}{c_2t}}.
 \end{equation}

Recall that the notion of atoms introduced in \cite{Preisner_Sikora_Yan} was different from the atoms introduced in Definition~\ref{def_atoms}. The atoms from \cite{Preisner_Sikora_Yan} are denoted 
as $[\mu,h]$-atoms in contrast to $(\mu;h;2)$-atoms that we consider here.
The $[\mu,h]$-atom is defined by the following condition
\eq{\label{old_atoms}
\supp \, a \subseteq B, \quad \norm{a}_{L^2\eee{h^{-1}\mu}} \leq \mu_{h}(B)^{-1/2}, \quad \int a(x)h(x)\, d\mu(x)=0,
}
where $\mu_h$ is the measure on $X$ with the density $h(x) d\mu(x)$. Notice that the difference between $[\mu,h]$-atoms and $(\mu;h;2)$-atoms lies in the $L^2$-size condition. The space $H^1_{at}[\mu,h]$ is defined analogously as in Definition \ref{def_atoms}. It appears that, under suitable assumptions, the spaces $H^1_L(X)$, $H^1_{at}(\mu;h;2)$, and $H^1_{at}[\mu,h]$ are equivalent. 

At this point let us also mention that the assumption \eqref{L-harm} is\blue{, in some sense, } superfluous. If a function $h$ satisfies \eqref{ULG'} then one can find $\wt{h}$ that is comparable to $h$ and satisfy \eqref{L-harm}, see \cite[Proposition 2.3]{Preisner_Sikora_Yan}.

In order to state our results, recall that that a non-negative function $w\in L^1_{loc}(\mu)$ on $X$ is in the Reverse H\"older class $RH_q(\mu)$, $1< q <\infty$, if there is a constant $C>0$ such
that
\begin{equation}\label{RHq}
    \left(\frac{1}{\mu(B)} \int_B w(x)^q d\mu(x) \right)^{1/q}
\leq \frac{C}{\mu(B)}\int_B w(x) d\mu(x)
\end{equation}
for every ball $B$ in $X$. Moreover, $RH_\infty(\mu)$ is defined by the condition
\eq{
\label{RHinfty} \sup_{x\in B} w(x) \leq \frac{C}{\mu(B)} \int_B w(x)\, d\mu(x).
}
{
Recall that  $RH_{q_2}(\mu) \subseteq RH_{q_1}(\mu)$ for $1<q_1\leq q_2 \leq \8$.
}

We are ready to state a theorem that establishes the equality $H^1_L(X) = H^1_{at}(\mu;h;2)$.

\sthm{main_thm_1}{
Suppose that a metric measure space $(X,d,\mu)$ satisfies the doubling condition \eqref{doubling}. Let $L$ be a self-adjoint densely defined operator on $L^2(X)$, such that the the upper Gaussian bounds \eqref{UG} holds for the corresponding heat kernel. Assume that: $h: X \to [0,\8)$ is an $L$-harmonic function in the sense of \eqref{L-harm}, the estimates \eqref{ULG'} are satisfied, and $h\in RH_\infty(\mu)$. Then the spaces $H^1_L(X)$ and $H^1_{at}(\mu;h;2)$ coincide and have comparable norms.
}

The proof of Theorem \ref{main_thm_1} will be given in Section \ref{sec2} below. We additionally prove a similar result with slightly weaker assumption, namely $h\in RH_q(\mu)$ with $q<\infty$ large enough, see Theorem~ \ref{main_thm_2}.
Let us point out the differences between Theorem \ref{main_thm_1} and \cite[Theorem A]{Preisner_Sikora_Yan}, which we recall in Theorem \ref{old_thm} below. Firstly and most significantly, the $(\mu;h;2)-$atoms in Theorem~\ref{main_thm_1} satisfy $\norm{a}_{L^2\eee{\mu}} \leq \mu(B)^{-1/2}$, whereas $[\mu,h]$-atoms from \cite{Preisner_Sikora_Yan} have a more complicated size condition $\norm{a}_{L^2\eee{h^{-1}\mu}} \leq \mu_h(B)^{-1/2}$.
The later involves the function $h$. We shall see that this change opens a possibility of studying manifolds in the case when there exist two or more linearly independent bounded  harmonic  functions. It also has clearer geometrical meaning when $h$ is used only in cancellation condition. Secondly, the assumption $h^{-1} \in A_p(\mu_{h^2})$ from Theorem \ref{old_thm} is replaced by $h\in RH_\8(\mu)$. Thirdly, a simple fact stated in Proposition~\ref{lem_doubling} allows us to get rid of the assumption that $\mu_{h^2}$ is a measure that satisfies the doubling property. Let us mention that for all the applications that are familiar to us the condition $h\in RH_\8(\mu)$ holds.

As a consequence of our new atomic characterization of $H^1_L(X)$ we  also obtain a more natural description of the dual space.
\sthm{thm_duality}{
Assume that $h\in L^2_{loc}(\mu)$. Then the space $BMO(\mu;h;2)$ is dual to $H^1_{at}(\mu;h;2)$.
}

Let us notice that the above statement does not refer to the operator $L$ and the corresponding semigroup $T_t$. However, its main motivation comes from Theorems \ref{main_thm_1}. Theorem \ref{thm_duality} will be proved in a more general setting (with a set of functions $h_1,...,h_m$ instead of just one function $h$) in Section \ref{sec3}.

\subsection{Two harmonic functions}\label{sec13}

Another goal of this paper is to study operators that are related to two (or more) linearly independent harmonic functions. We shall discuss this situation in a particular case of manifolds with ends, but one can apply our result in other similar circumstances. Analysis on manifolds with ends has recently attracted a significant board interest. The theory of heat kernel bounds on such manifolds was studied, among others, by Davies, Grigor'yan, Ishiwata, and Saloff-Coste in \cite{Davies2, Grigoryan_Saloff-Coste, G_I_S-C1, G_I_S-C2}. The Riesz transform in this context was studied in, e.g. \cite{Sikora_Bailey, Dangyang, Carron_Coulhon_Hassell,
Sikora_Hassell, Sikora_Nix}.

We shall consider  a symmetric Riemannian  manifold with two ends of the form $M = \Rn \# \Rn$, $n\geq 3$. Such manifold is obtained as a~connected sum of two copies of $\Rn$ combined by a~compact smooth connection.
 We shall also assume that there exists a symmetry $\sigma: M \to M$ that interchanges one end with the another and preserves both: the metric and the measure. For precise descriptions, definitions, and assumptions see Section \ref{sec5.1}.

Let $L_M$ be a Laplace-Beltrami operator on $M = \Rn \#\Rn$. Denote the associated heat semigroup by $T_t = \exp(-tL_M)$. It is known that the corresponding heat kernel  $T_t(x,y)$  satisfies \eqref{UG}, see e.g. \cite{Grigoryan_Saloff-Coste}. In addition, for $L_M$ there exist two linearly independent non-negative bounded harmonic functions $h_+$ and $h_-$. One possible choice of $h_+$ can be given in probability language, where $h_+(x)$ is the probability that the process starting at $x\in M$ will escape at $t\to \8$ via ''the right end''. For analytic description of these functions, see e.g. \cite{Sikora_Hassell}.

The following theorem establishes the equivalence of $H^1_{L_M}(M)$ and the atomic Hardy space $H^1_{at}(\mu; h_+,h_-;2)$, see Definition~\ref{def_atoms}.

\sthm{thm_manifold}{
Let $L_M$ be the Laplace-Beltrami operator on a symmetric manifold $M=\Rn\#\Rn$, see Section \ref{sec5.1}. Denote by $h_+$ and $h_-$ the bounded $L_M$-harmonic functions. Then the Hardy spaces $H^1_{L_M}(M)$ and $H^1_{at}(\mu; h_+, h_-; 2)$ coincide and there exists $C>0$ such that
$$C^{-1} \norm{f}_{H^1_{L_M}(M)} \leq \norm{f}_{H^1_{at}(\mu; h_+,h_-;2)}\leq C\norm{f}_{H^1_{L_M}(M)}.$$
}

The proof of Theorem \ref{thm_manifold} is discussed in Section \ref{sec5.1} below. Using Theorem~\ref{thm_duality_general} one can describe the dual space as follows.

\cor{cor_duality}{
The dual space of the atomic Hardy space $H^1_{at}(\mu;h_+,h_-;2)$ from Theorem \ref{thm_manifold} is the  space $BMO(\mu;h_+,h_-;2)$.
}

\section{Atomic characterizations of $H^1_L(X)$.}\label{sec2}

In this section we consider the metric measure space $(X,d,\mu)$ and operators $L, T_t$ as in Sections \ref{ssec11} and \ref{ssec12}. In particular, we assume that there exists $h$ such that \eqref{L-harm} and \eqref{ULG'} are satisfied. Our main goal here is to prove Theorem \ref{main_thm_1}. Actually, we shall prove a slightly more general Theorem \ref{main_thm_2} that covers Theorem \ref{main_thm_1} as a special case $q=\8$ and $r=2$.

\thm{main_thm_2}{
Suppose that a metric measure space $(X,d,\mu)$ satisfies the doubling condition \eqref{doubling}. Let $L$ be a self-adjoint operator on $L^2(X)$, such that the the upper Gaussian bounds \eqref{UG} holds for the corresponding heat kernel semigroup. Assume that: $h: X \to [0,\8)$ is an $L$-harmonic function in the sense of \eqref{L-harm}, the estimates \eqref{ULG'} are satisfied, and $h\in RH_q(\mu)$ for $q\in (3,\8]$ large enough. Then we  have that $H^1_L(X)=H^1_{at}(\mu;h;r)$ with $r=2q/(q+1)$ if $q<\8$ and $r=2$ if $q=\8$. Moreover the corresponding norms are equivalent.
}

In this note the parameters: $q\in(1,\8]$ from $RH_q(\mu)$, see \eqref{RHq}, and $r\leq 2$ from $(\mu;h;r)$-atoms are always related by the equation
\eq{\label{eq_rq}
r = \frac{2q}{q+1}.
}
Notice that if $r$ and $q>3$ satisfy  \eqref{eq_rq}, then $r'<q$. Thus $h\in RH_q(\mu)$ implies that $h\in RH_{r'}(\mu)$, and $h\in L^{r'}(\mu)$ as it is required in Definition \ref{def_atoms}. 
There is a natural question whether the spaces $H^1_{at}(\mu;h;r)$ are equivalent for a full range $r\in(1,\8]$, but we do not investigate this issue here.

\subsection{Auxiliary results}

Before going to the proofs of Theorems \ref{main_thm_1} and \ref{main_thm_2} we shall recall some results and prove auxiliary estimates.

By definition, a locally integrable non-negative function $w$ is in the Muckenhoupt class $A_p(\nu)$ if 
\eq{
\label{muck}
\eee{\frac{1}{\nu(B)} \int_B w(x) d\nu(x)  }    \eee{\frac{1}{ \nu(B)} \int_B w(x)^{-1/(p-1)}(x)d\nu(x)}^{p-1} \leq C
}
for all balls $B$.

Recall that the atomic Hardy  space $H^1_{at}[\mu,h]$ with $[\mu,h]$-atoms (see \eqref{old_atoms}) are considered in \cite{Preisner_Sikora_Yan}, where the following theorem is proved.

\thm{old_thm}{\cite[Theorem A]{Preisner_Sikora_Yan}
Let $(X,\mu,d)$ be a space satisfying \eqref{doubling}, $L$ be a non-negative self-adjoint operator related to the semigroup $T_t$ with a kernel $T_t(x,y)$ satisfying \eqref{UG}. We assume that there exists a function $h:X \to (0,\8)$ such that: $h$ is $L$-harmonic in the sense of \eqref{L-harm}, the measure $\mu_{h^2}$ is doubling, and \eqref{ULG'} is satisfied. Then there exists $p_0\in (1,2]$ such that if $h^{-1}\in A_{p_0}(\mu_{h^2})$, then $H^1_L(X)$ and $H^1_{at}[\mu,h]$ coincide and
$$\norm{f}_{H^1_L(X)} \simeq \norm{f}_{H^1_{at}[\mu,h]}.$$
}

Recall the notion of so-called Doob transform (or $h$-transform) technique, see e.g. \cite[Section 5.1.2]{Gyrya_Saloff-Coste}. This is one of the main tools in the proof of Theorem \ref{old_thm}. For the semigroup $T_t$ on $L^2(\mu)$ one can define a new kernel
$$\wt{T}_t(x,y) = \frac{T_t(x,y)}{h(x)h(y)},$$
which defines a semigroup on the space $L^2(\mu_{h^2})$. The semigroup $\wt{T}_t$ is conservative, i.e. $\int_X \wt{T}_t(x,y) \, d\mu_{h^2}(x)=1$. Notice that the assumption \eqref{ULG'} translates into the upper and lower Gaussian bounds for $\wt{T}_t(x,y)$, i.e.
\eq{\label{eq1234}
C^{-1}\mu_{h^2}(B(x,\sqrt{t}))^{-1} \exp\eee{-\frac{d(x,y)^2}{c_1t}} \leq \wt{T}_t (x,y) \leq C\mu_{h^2}(B(x,\sqrt{t}))^{-1} \exp\eee{-\frac{d(x,y)^2}{c_2t}}.
}
It is known that the estimates \eqref{ULG'} automatically give the following H\"older-type  estimates,  for a short proof see e.g. \cite[Section 4]{DP_Annali}.
See  also \cite[Theorem 5.11]{Gyrya_Saloff-Coste}. 
\prop{lem_Holder}{
Assume that $(X,d,\mu)$ is a metric measure space that satisfies \eqref{doubling}. Let $T_t$ be a semigroup related to a kernel $T_t(x,y)$ such that \eqref{ULG'} holds for a function $h$. Then there exist $C, \delta >0$ such that for $t>0$ and $x,y,y_0\in X$ satisfying $d(y,y_0)\leq \sqrt{t}$ we have
\eq{\label{Holder}
\abs{\frac{T_t(x,y)}{h(x)h(y)} - \frac{T_t(x,y_0)}{h(x)h(y_0)}} \leq C \eee{\frac{d(y,y_0)}{\sqrt{t}}}^\delta \mu_{h^2}(B(x,\sqrt{t}))^{-1} \exp\eee{-\frac{d(x,y)^2}{c_3t}}. 
}
}
\begin{proof}
Note that $\int_X \wt{T}_t(x,y)\, d\mu_{h^2}(x)=1$ a.e.  and the semigroup $\wt{T}_t$ satisfies \eqref{eq1234}. Therefore \eqref{Holder} follows from \cite[Theorem 5]{DP_Annali}.
\end{proof}

Another known consequence of \eqref{ULG'} is the following proposition.
\prop{lem_doubling}{
Assume that a semigroup $T_t$ related to a kernel $T_t(x,y)$ satisfies \eqref{ULG'} for some function $h$. Then the measure $\mu_{h^2}$  satisfies the doubling condition \eqref{doubling}.
}
\pr{
The lemma is well-known and follows from \eqref{eq1234}. We shall provide a short proof for the sake of completeness. Let $x\in X$, $r>0$. By \eqref{eq1234}  for $y\in B(x,r)$ we have that
$\mu_{h^2}(B(x,\sqrt{2}r))^{-1} \gtrsim \wt{T}_{2r^2}(x,y)$. By averaging this inequality over $B(x,r)$ and using the semigroup property we obtain
\spx{
\mu_{h^2}(B(x,\sqrt{2}r))^{-1} &\gtrsim   \mu_{h^2}(B(x,r))^{-1} \int_{B(x,r)} \wt{T}_{2r^2}(x,y)\, d\mu_{h^2}(y) \\
&\gtrsim  \mu_{h^2}(B(x,r))^{-1} \int_{B(x,r)} \int_{B(y,r)} \wt{T}_{r^2}(x,z) \wt{T}_{r^2}(z,y)\,d\mu_{h^2}(z) \, d\mu_{h^2}(y)\\
&\gtrsim  \mu_{h^2}(B(x,r))^{-2} \int_{B(x,r)} \int_{B(y,r)} \mu_{h^2}(B(y,r))^{-1}\,d\mu_{h^2}(z) \, d\mu_{h^2}(y)\\
&\simeq  \mu_{h^2}(B(x,r))^{-1},
}
which ends the proof.
}

The next lemma describes relation between $[\mu,h]$-atoms and $(\mu;h;r)$-atoms. This is required in the first part of the proof of Theorems \ref{main_thm_1} and \ref{main_thm_2}.
\lem{lem_atoms}{
{\it (a)} Assume that $h\in RH_\8(\mu)$. Then there exists $C$  such  that for a ball $B$ and a  function $a$ supported in $B$ we have
\spx{
\norm{a}_{L^2\eee{h^{-1}\mu}} \leq \mu_h(B)^{-1/2} \implies \norm{a}_{L^2\eee{\mu}} \leq C \mu(B)^{-1/2}.
}
{\it (b)} Assume that $h\in RH_q(\mu)$ for $q>1$ and $r<2$ is as in \eqref{eq_rq}. Then there exists $C>0$  such that for any ball $B$ and any  function $a$ supported in $B$ we have
\spx{
\norm{a}_{L^2\eee{h^{-1}\mu}} \leq \mu_h(B)^{-1/2} \implies \norm{a}_{L^r\eee{\mu}} \leq C \mu(B)^{-1+1/r}.
}
}
\pr{
The part {\it (a)} follows immediately  from the definition of $RH_\8(\mu)$. Indeed
\spx{
\norm{a}_{L^2(\mu)}^2 \leq \sup_{x\in B} h(x) \cdot \norm{a}_{L^2(h^{-1}\mu)}^2 \leq \sup_{x\in B} h(x) \cdot \mu_h(B)^{-1} \lesssim \mu(B)^{-1}.
}

To prove {\it (b)} we use H\"older's inequality with the exponents $2/r$ and $2/(2-r)$ and the definition of $RH_q(\mu)$ with $q=r/(2-r) >1$, see \eqref{RHq}. Then
\spx{
\norm{a}_{L^r(\mu)}^r &= \int_B |a|^r h^{-r/2} h^{r/2} \, d\mu \leq \eee{\int_B |a|^2 h^{-1} \, d\mu}^{r/2} \cdot \eee{\int_B h^{r/(2-r)}\,d\mu}^{(2-r)/2}\\
&\lesssim \mu_h(B)^{-r/2} \cdot \mu(B)^{1-r} \cdot \mu_h(B)^{r/2} =\mu(B)^{1-r}.
}
}

In the whole paper $p$ and $q$ denote the parameters from the Muckenhoupt class $A_p$ and the reverse H\"older class $RH_q$, respectively. They are always related by the equation
\eq{\label{eq_pq}
p=\frac{q-1}{q-2}.
}

\lem{lem_Muc-RH}{
Assume that $h\in RH_q(\mu)$ for $q>2$. Then we have $h^{-1} \in A_p(\mu_{h^2})$ for $p$ as in \eqref{eq_pq}. In particular, if $h\in RH_\8(\mu)$, then $h^{-1} \in A_p(\mu_{h^2})$ for all $p>1$.
}
\pr{
Let $p>1$ and $q>2$ satisfy \eqref{eq_pq}, i.e. $q=2+(p-1)^{-1}$. Then the condition $h^{-1}\in A_p(\mu_{h^2})$ is equivalent to
\eq{\label{pii}
\mu(B)^{-1} \int_B h \, d\mu \cdot \eee{\mu(B)^{-1} \int_B h^q d\mu}^{p-1}\lesssim \eee{\mu(B)^{-1} \int_B h^2 \, d\mu}^{p},
}
see \eqref{muck}. By the Reverse H\"older's assumption the left-hand side of \eqref{pii} can be estimated by $ \eee{\mu(B)^{-1}\int_B h \, d\mu}^{2p}$. Then, \eqref{pii} follows from the Cauchy-Schwarz inequality.
 The second statement follows immediately, since if $q\to \8$ then $p=1+1/(q-2)$ tends to $1$.
}

\subsection{Proofs of Theorems \ref{main_thm_1} and \ref{main_thm_2}}

\pr{[Proof of Theorem \ref{main_thm_2}]
{\bf Case 1: $q<\8$.} Let $q_0 = \max(D/\delta, 2+(p_0-1)^{-1})$, where $D$ is from \eqref{doubling}, $\delta$ is the H\"older exponent from \eqref{Holder}, and $p_0$ is as in Theorem \ref{old_thm}. Notice that since $p_0\leq 2$ then $q_0 \geq 3$. Assume that $h\in RH_q(\mu)$ with $q>q_0$ and, as before, $r=2 -2/(q+1)$, see \eqref{eq_rq}.

{\bf Step 1.} Since $q>q_0$, then $h\in RH_{q_0}(\mu)$ and, by Lemma \ref{lem_Muc-RH}, we obtain that $h^{-1}\in A_{p_0}(\mu_{h^2})$. Recall that, by Proposition \ref{lem_doubling}, the measure $\mu_{h^2}$ is doubling. By Theorem \ref{old_thm}, we have that $H_L^1(X) = H^1_{at}[\mu,h]$. Moreover, Lemma \ref{lem_atoms} {\it (b)} states that every $[\mu,h]-$atom is also an $(\mu,h,r)$-atom, so $H^1_L(X) \subseteq H^1_{at}(\mu,h,r)$.

{\bf Step 2a.} In this step we prove that there exists $C>0$ such that
\eq{\label{atom_bdd}
\norm{M_La}_{L^1(\mu)} \leq C
}
for every $(\mu;h; r)$-atom $a$ with $C$ independent of $a$. Assume that $a$ is related to a ball $B=B(y_0,R)$ and satisfies \eqref{new_atoms}. Denote $\lambda B = B(y_0, \lambda R)$, $\lambda>0$. By \eqref{doubling} and the boundedness of $M_L$ on $L^r(\mu)$  we have
$$\norm{M_La}_{L^1(2B,\, \mu)} \leq \mu(2B)^{1-\frac{1}{r}} \norm{M_La}_{L^r(\mu)}\lesssim \mu(B)^{1-\frac{1}{r}} \norm{a}_{L^r(\mu)} \leq C.
$$

Our next goal is to estimate $M_La(x)$ for $x\not\in 2B$. Notice that if $y\in B$, then $d(x,y) \simeq d(x,y_0)$. Consider first $t$ such that \underline{$\sqrt{t}>d(y,y_0)$}. In this case, the cancellation $\int a h \, d\mu =0$ and \eqref{Holder} yields
\sp{\label{eq45678}
|T_t a(x)| &\leq h(x) \int_B \abs{\frac{T_t(x,y)}{h(x)h(y)}-\frac{T_t(x,y_0)}{h(x)h(y_0)}}|a(y)| h(y)\, d\mu(y)\\
&\lesssim  h(x) \int_B \eee{\frac{d(y,y_0)}{\sqrt{t}}}^{\delta} \mu_{h^2}(B(x,\sqrt{t}))^{-1} \exp\eee{-\frac{d(x,y)^2}{ct}}|a(y)|h(y) \, d\mu(y)\\
&\lesssim h(x) \eee{\frac{R}{\st}}^\delta  \mu_{h^2}(B(x,\sqrt{t}))^{-1} \exp\eee{-\frac{d(x,y_0)^2}{ct}}\int_B |a| h \, d\mu\\
&\lesssim h(x) \eee{\frac{R}{\st}}^\delta  \mu_{h^2}(B(x,d(x,y_0)))^{-1} \exp\eee{-\frac{d(x,y_0)^2}{c't}}\int_B |a| h \, d\mu,
}
Next suppose  \underline{$\sqrt{t}\leq d(y,y_0)$}, we have that $R>\sqrt{t}$ and, by \eqref{ULG'}, we get
\sp{\label{eq5678}
|T_t a(x)| &\leq h(x) \int_B \frac{T_t(x,y)}{h(x)h(y)}|a(y)| h(y)\, d\mu(y)\\
&\lesssim h(x) \eee{\frac{R}{\st}}^\delta \mu_{h^2}(B(x,\sqrt{t}))^{-1} \exp\eee{-\frac{d(x,y_0)^2}{ct}} \int_B  |a|h \, d\mu.
}
Notice that Proposition \ref{lem_doubling} implies
\sp{\label{eq56789}
 \frac{ \mu_{h^2}(B(x,d(x,y_0)))}{ \mu_{h^2}(B(x,\sqrt{t}))}  \lesssim \eee{1+\frac{d(x,y_0)}{\sqrt{t}}}^D.
}
By combining \eqref{eq45678}, \eqref{eq5678}, \eqref{eq56789}, and using  H\"older's inequality we arrive at
\sp{\label{eq22333}
M_La(x) &\lesssim  h(x) \sup_{t>0}\eee{ \eee{\frac{R}{\st}}^\delta  \mu_{h^2}(B(x,\sqrt{t}))^{-1} \exp\eee{-\frac{d(x,y_0)^2}{ct}}}\int_B |a| h \, d\mu\\
&\lesssim h(x) \mu_{h^2}(B(x,d(x,y_0)))^{-1}  \sup_{t>0}\eee{\eee{\frac{R}{\st}}^\delta   \exp\eee{-\frac{d(x,y_0)^2}{c't}} }\int_B |a| h \, d\mu,\\ &\lesssim \eee{\frac{R}{d(x,y_0)}}^\delta \frac{h(x)}{\mu_{h^2}(B(x,d(x,y_0)))} \norm{a}_{L^{r}(B,\, \mu)} \norm{h}_{L^{r'}(B,\, \mu)},
}
where $1/r + 1/r'=1$.

For $n=1,2,\ldots$ denote $S_n = 2^{n+1}B \setminus 2^nB$. Then
$$\norm{M_La}_{L^1((2B)^c,\, \mu)}\leq \sum_{n=1}^\8 \underbrace{\norm{M_La}_{L^1(S_n,\, \mu)}}_{P_n} $$
We claim that $\sum_{n\geq 1} P_n \leq C$. Notice that $\mu_{h^2}(B(x,d(x,y_0)))\simeq \mu_{h^2}(B(y_0,d(x,y_0)))$, see Proposition \ref{lem_doubling}. For $x\in S_n$ we have that $2^nR \leq d(x,y_0) \leq 2^{n+1}R$, so by \eqref{eq22333} and \eqref{new_atoms} we obtain
\spx{
P_n &\lesssim 2^{-n\delta}  \frac{\mu_h(2^{n+1}B)}{  \mu_{h^2}(2^nB)}
\mu(B)^{-1+1/r}\norm{h}_{L^{r'}(B,\, \mu)} \\
&\simeq 2^{-n\delta} \frac{\eee{\mu(2^{n+1}B)^{-1}\int_{2^{n+1}B} h\, d\mu} \cdot \eee{\mu(B)^{-1}\int_{B} h^{r'}\, d\mu}^{1/r'}}{\mu(2^{n+1}B)^{-1}\int_{2^{n+1}B} h^2\, d\mu}\\
&\simeq 2^{-n\delta} \frac{ \eee{\mu(B)^{-1}\int_{B} h^{r'}\, d\mu}^{1/r'}}{\eee{\mu(2^{n+1}B)^{-1}\int_{2^{n+1}B} h^2\, d\mu}^{1/2}}.
}
In the estimate above we used Proposition \ref{lem_doubling} and H\"older's inequality. Recall that $q>q_0\geq 3$, $r'=2q/(q-1)\leq q$, and $h\in RH_q(\mu)$. By H\"older's inequality and Reverse H\"older's assumption we arrive at
\spx{
P_n&\lesssim 2^{-n\delta}  \frac{ \eee{ \mu(B)^{-1}\int_{B} h^{q}\, d\mu}^{1/q}}{\mu(2^{n+1}B)^{-1} \int_{2^{n+1}B} h\, d\mu}\\
&\lesssim 2^{-n\delta}  \eee{\frac{\mu(2^{n+1}B)}{\mu(B)}\cdot \frac{\int_B h^q\, d\mu}{\int_{2^{n+1}B}h^q\, d\mu}}^{1/q}\\
&\lesssim 2^{-n\delta} 2^{nD/q},
}
where $D$ is from \eqref{doubling}. The claim $\sum_{n\geq 1} P_n \leq  C$ is proved, since $\delta -D/q>0$. The proof of \eqref{atom_bdd} is finished.

\label{page7}
{\bf Step 2b.} The inclusion $H^1_{at}(\mu;h;r) \subseteq H_L^1(X)$ follows from \eqref{atom_bdd} by a standard argument. We present the details for the sake of completeness. Take $f(x)=\sum_k \la_k a_k(x)$, where $\sum_k |\la_k|<\8$ and $a_k$ are $(\mu;h;r)$-atoms. Recall that $H^1_L(X)$ embeds continuously into $L^1(\mu)$, see \cite{AMM} and \cite[Section 8]{Preisner_Sikora_Yan}. Let $f_n(x) = \sum_{k=1}^n \la_k a_k(x)$. By \eqref{atom_bdd} for $n<m$ we have
$$\norm{M_L(f_m-f_n)}_{L^1(\mu)} \lesssim \sum_{k=n+1}^m |\la_k|.$$
Thus $f_n$ is a Cauchy sequence in $H^1_L(X)$ and there is a limit $g\in H^1_L(X)$ in the $H^1_L(X)$-norm. Also, $f_n \to f$ in $L^1(\mu)$ and, since $H^1_L(X)$ embeds continuously into $ L^1(\mu)$, we have that $f=g$. Moreover,
$$\norm{f}_{H^1_L(X)} \leq \sup_n \norm{f_n}_{H^1_L(X)} \lesssim \sum_{k=1}^\8 |\la_k| \lesssim \norm{f}_{H^1_{at}(\mu;h;2)},$$
which ends the proof of Theorem \ref{main_thm_2} in the case $q<\8$.
}

{\bf Case 2: $q=\8$.} This case (equivalent with Theorem \ref{main_thm_1}) can be verified by the same argument with appropriate modifications and we only comment on the differences. The assumption $h\in RH_\8(\mu)$ implies $h\in RH_q(\mu)$ for all $q>1$. Then $h^{-1}\in A_p(\mu_{h^2})$ for all $p>1$, see Lemma~ \ref{lem_Muc-RH}. The main change in the proof is that we  use Lemma \ref{lem_atoms} {\it (a)} (instead of {\it (b)}), and, therefore, we may keep $r=2$. We leave the details to the reader.

\section{BMO spaces and harmonic functions.}\label{sec3}

Let $r,r'>1$ be the conjugate exponents, i.e. $1/r+1/r' = 1$.
In this section we prove that the dual of the atomic Hardy space $H^1_{at}(\mu;h_1,...,h_m;r)$, see Definition \ref{def_atoms},  is the space $BMO(\mu;h_1,...,h_m;r')$, see Definition \ref{BMO_general}.

\thm{thm_duality_general}{
Let $m\in \NN$ and $h_1,...,h_m\in L^{r'}_{loc}(\mu)$ be functions that are  linearly independent on $X$. Then the space $BMO(\mu;h_1,...,h_m;r')$ is dual to $H^1_{at}(\mu;h_1,...,h_m;r)$.
}

The proof of Theorem \ref{thm_duality_general} consists of two statements:
\en{
\item[{\bf 1)}] For each $g\in \BMO$ one can define (in a linear way) a bounded linear functional on $H^1_{at}(\mu;h_1,...,h_m;r)$ that we shall denote by $\langle f,g \rangle$, such that
$$|\langle f,g \rangle | \leq \norm{f}_{H^1_{at}(\mu;h_1,...,h_m;r)} \norm{g}_{\BMO}.$$
\item[{\bf 2)}] Conversely, for every $l \in (\Har)^*$ there exists $g\in \BMO$, such that $l$ is defined as above, i.e. $l(f) = \langle f, g \rangle$ with $$\norm{g}_{BMO(\mu;h_1,...,h_m;r)} \lesssim \norm{l}_{H^1_{at}(\mu;h_1,...,h_m;r) \to \CC}.$$
}

\pr{
Theorem \ref{thm_duality} follows by a proper modification of standard techniques, see e.g. \cite[p.~142-145]{Stein}, \cite[Section 5]{Preisner_Sikora_Yan}. Denote $V = \mathrm{span}(h_1,...,h_m)$ as in Definition \ref{BMO_general}.

To prove {\bf 1)} take $g\in \BMO$. For an $(\mu;h_1,...,h_m;r)$-atom $a$ and $u\in V$ we have that
\spx{
\abs{\int_X a(x) g(x) d\mu(x)} &= \abs{\int_X a(x) (g(x)-u(x)) d\mu(x)}\\
&\leq \norm{a}_{L^r(\mu)} \eee{\int_B \abs{g(x) -u(x)}^{r'}  d\mu(x)}^{1/{r'}}\\
&\leq \eee{\mu(B)^{-1}\int_B \abs{g(x) -u(x)}^{r'}  d\mu(x)}^{1/{r'}}.
}
The estimate above holds for all balls $B$ in $X$ and all $u\in V$, therefore we have that
\eq{\label{est1222}
\abs{\int_X a(x) g(x) d\mu(x)} \leq \norm{g}_{\BMO}.
}
Let $f(x)=\sum_k \la_k a_k(x)$, where $\sum_k|\la_k|<\8$. From \eqref{est1222} we deduce that
$$\int_X \sum_{k=1}^n \la_k a_k(x) g(x), d\mu(x)$$
is a Cauchy sequence as $n\to \8$ and we define
the pairing $\langle f, g \rangle$ as its limit for an arbitrary $f\in \Har$. It is easy to check that this pairing is well-defined, linear, and 
$$|\langle f,g\rangle| \leq \norm{f}_{\Har} \norm{g}_{\BMO}.$$
Notice that $\langle f,g \rangle = \int f (x) g(x) \, d\mu(x)$ when $f$ is a finite linear combination of $(\mu;h_1,...,h_m;r)$-atoms (this subspace is dense in $\Har$).

Now we turn to prove {\bf 2)}. Let $l$ be a linear functional on $\Har$.We may assume that $\norm{l}_{\Har \to \CC} \leq 1$.
Let $B$ be a ball and consider the space
$$\hh_B = \set{f \in L^r(\mu|_B) \, : \, \forall u\in V\ \   \int_B f(x) u(x) d\mu(x) = 0}.$$
Obviously, if $f\in \hh_B$, then $f\in  \Har$ with
$$\norm{f}_{\Har} \leq \mu(B)^{1-1/r} \norm{f}_{L^r(\mu)}.$$
Therefore, by the Hahn-Banach extension and the Riesz representation theorems, there exists $\wt{g}_B$ such that
$$l(f) = \int_B f(x)  \wt{g}_B(x)  \, d\mu(x), \quad f\in \hh_B,$$
$$\norm{\wt{g}_B}_{L^{r'}(\mu)}=\norm{l}_{\hh_B \to \CC} \leq \mu(B)^{1-1/r}.$$

Since the functions $h_1,...,h_m$ are linearly independent on $X$, then there exists a ball such that $h_1,...,h_m$ are linearly independent on this ball. Fix such a ball $B_0$. Then $M=(\int_{B_0} h_i h_j \, d\mu)_{i,j=1}^m$ is a non-degenerate $m\times m$ matrix. We choose the constants $c_B = (c_{B,1},...,c_{B,m})$ so that
$$M c_B = -\eee{\int_{B_0} \wt{g}_B \, h_1\, d\mu, ..., \int_{B_0} \wt{g}_B \, h_m\, d\mu}.$$
We set $g_B = \wt{g}_B + c_{B,1}h_1+...+c_{B,m} h_m$. Notice that $l(f)=\int_B f(x) g_B(x)\, d\mu(x)$ for $f\in \hh_B$. According to our choice of constants $c_B$ we have
\eq{\label{skr123}
\int_{B_0} g_B \, h_j\, d\mu=0, \qquad j=1,...,m.
}
Take an increasing family  of balls $B_0 \subseteq B_1 \subseteq ...$ {such that $\bigcup_{n\in \NN} B_n = X$}.
We have that $g_{B_n} - g_{B_{n+1}} = d_{n,1} h_1 + ...+d_{n,m} h_m$ on $B_n \supseteq B_0$. By integrating on $B_0$ with respect to $h_1,...,h_n$ we obtain that $M (d_{n,1},...,d_{n,m})=0$ and thus $d_{n,j} =0$ for $n\in \NN$ and $j=1,...,m$. Since $g_{B_n}$ and $g_{B_{n+1}}$ coincide on $B_n$ we can define the limit
$$g(x) =  \lim_{n\to \infty} g_{{B_n}}(x).$$
Moreover, for any ball $B$,
\begin{equation}
\eee{\frac{1}{\mu(B)} \int_B \abs{g(x)-c_{B,1} h_1 -...-c_{B,m} h_m }^{r'} d\mu(x)}^{\frac{1}{r'}} 
= \eee{\frac{1}{\mu(B)} \norm{\wt{g}_B}_{L^{r'}(\mu)}^{r'} }^{\frac{1}{r'}} \leq C.
\end{equation}
We have proved  that $g\in \BMO$ and $\norm{g}_{\BMO} \leq C$. Also, $l(f) = \int_X f(x) g(x) d\mu(x)$ whenever $f$ is a finite
 combination of atoms. This ends the proof of {\bf 2)}.
}

\section{Spaces with symmetry and operators related to two harmonic functions}
\label{sec4}

Let $(X,d,\mu)$ be a space of homogeneous type, see \eqref{doubling}. Let us start with a general proposition, that relies on \eqref{UG}. Additional assumptions (e.g. symmetry of $X$) will be added later.

\prop{prop_skr}{
Suppose that  $L$ is a non-negative self-adjoint operator acting  on $L^2(\mu)$ such that the integral kernel  $T_t(x,y)$ corresponding to 
  the semigroup $T_t=e^{-tL}$  satisfies the upper Gaussian estimates \eqref{UG}. Assume that the functions $h_1,...,h_m \in L^\8_{loc}(\mu)$ are $L$-harmonic in the sense of \eqref{L-harm}. Then $H^1_L(X) \subseteq H^1_{at}(\mu;h_1,...,h_m;2)$ and there exists $C>0$, such that
$$\norm{f}_{H^1_{at}(\mu;h_1,...,h_m;2)} \leq C \norm{f}_{H^1_L(X)}.$$
}

\begin{proof}
The proof follows from  \cite{Hofmann_Memoirs} and the argument is similar to \cite[Lemma 9.1]{Hofmann_Memoirs} and \cite[Theorem 3.5]{Preisner_Sikora_Yan}. For the completeness we present a sketch of the proof.

Recall that the space $H^1_L(X)$ with $T_t$ satisfying \eqref{UG} has many equivalent characterizations, e.g. by maximal functions, Lusin area functions, square functions, etc. See \cite{Hofmann_Memoirs, Song_Yan_2018}. By \cite[Theorem 7.1]{Hofmann_Memoirs} for $f$ from a certain dense subset $\HH^1_L(X)$ of $H^1_L(X)$ we have that $f(x) = \sum_k \la_k(x) a_k$ (with the convergence in both $L^1(\mu)$ and $L^2(\mu)$), where $a_k$ are such that there exist balls $B_k = B(y_k,r_k)$ and functions $b_k \in \mathrm{Dom}  (L)$ such that:
$$a_k =Lb_k, \quad \supp \, (L^i b_k)\subseteq B_k, \quad \norm{(r_k^2L)^i b_k}_{L^2(\mu)}\leq r_k^2 \mu(B_k)^{-1/2}
$$
for $i=0,1$ and $k=1,2,...$. We can assume that $\norm{f}_{H^1_L(X)} \simeq \sum_k |\la_k|$. Since $a_k$ are properly localized and satisfy the size condition that we need, it is enough to show that
$\int a_k h_j\, d\mu  = 0$ for $j=1,...,m$ and $k\in \NN$.
For fixed $j$ from \eqref{L-harm} we have that
$$(I+L)^{-1} h_j = \int_0^\8 e^{-t} T_t h_j \, dt = h_j.$$
Then,
\spx{
\int_B a(x) h_j(x) \, d\nu(x) &= \int_B a(x) (I+\ll)^{-1}h_j(x) \, d\nu(x)\\
&= \int_B (I+\ll)^{-1}\ll b (x) h_j(x) \, d\nu(x)\\
&= \int_B (I+\ll)^{-1}(I+\ll)b (x) h_j(x) \, d\nu(x) - \int_B (I+\ll)^{-1}b (x) h_j(x) \, d\nu(x)\\
&= \int_B b (x) h_j(x) \, d\nu(x) - \int_B b (x) (I+\ll)^{-1} h_j(x) \, d\nu(x)\\
& = 0.
}
Notice that the calculation above is justified since $h_j\in L^\8_{loc}(\mu)$. What we have just proved is that each $a_k$ is $(\mu;h_1,...,h_m;2)$-atom. Thus, for $f\in \HH^1_L(X)$ we have that $f\in H^1_{at}(\mu;h_1,...,h_m;2)$ and
$$\norm{f}_{H^1_{at}(\mu;h_1,...,h_m;2)}\leq C \norm{f}_{H^1_L(X)}.$$
The proof is finished, since $\HH^1_L(X)$ is a dense subset of $H^1_L(X)$.
\end{proof}

\cor{coro_skr}{
Let $L$ be an operator as in Proposition \ref{prop_skr}. If $h$ is a bounded $L$-harmonic function then
$$\int f(x) h(x) \, d\mu(x) =0$$
for $f\in H^1_L(X).$
}
\begin{proof}
    Let $f\in H^1_L(X)$. From Proposition \ref{prop_skr} we have that $f(x) = \sum_k \la_k a_k(x)$, where $\sum_k |\la_k|<\8$ and $a_k$ are $(\mu;h;2)$-atoms. Then
    $$\int f(x) h(x) \, d\mu(x) = \lim_{n\to \8} \int \sum_{k=1}^n \la_k a_k(x) h(x)\, d\mu(x) = 0,$$
    since $\sum_{k=1}^n \la_k a_k(x)$ converges to $f$ in $L^1(\mu)$ and $h\in L^\8(\mu)$.
\end{proof}

\subsection{Background and assumptions}\label{ssec41}

\newcommand{\KKKK}{\Gamma}
In this section we shall consider spaces of homogeneous type with an additional symmetry. These assumptions are modelled by symmetric manifolds $\RR^n \# \RR^n$ with ends that are described in  Sections \ref{sec13} above and \ref{sec5.1} below. However the proofs do not rely on the manifold structure, thus we state some general assumptions that guarantee the atomic characterization of the Hardy space $H^1_L(X)$, where atoms have cancellations with respect to two linearly independent harmonic functions. According to our knowledge, such atomic Hardy spaces have not been intensively studied. However, see \cite{Paluszynski_Zienkiewicz}.

Let $(X,d,\mu)$ be a space of homogeneous type and $X = X_+ \cup X_-$.  Assume that $\KKKK:=X_+\cap X_-$, $\mu(\KKKK)=0$, and $\KKKK=\mathrm{bd} \, X_+=\mathrm{bd} \, X_-$. Denote by $\mu_+$ the measure $\mu$ restricted to $X_+$. We additionally assume that there exists an bijection $\sigma: X \to X$, that preserves both: the metric and the measure, $\sigma^2 = \mathrm{id}$, $\sigma(x) = x$ for $x\in \KKKK$, and $\sigma(X_+) = X_-$. In what follows we shall denote $\sigma(x) = \wt{x}$. We shall work under assumption that
\begin{equation}\label{metric}
\forall x,y\in X_+ \quad    d(x,y)\leq d(x,\wt{y}).
\end{equation}

We start by considering two operators on $X_+$ that satisfy the assumptions of Theorem~\ref{main_thm_1}. For further references let us recall these assumptions. Let $L_1$ and $L_2$  be two self-adjoint, non-negative operators that are densely defined on $L^2(\mu_+)$. Denote by $T_{t}^{[1]}$ and $T_{t}^{[2]}$ the semigroups related to $L_1$ and $L_2$, and let $T_{t}^{[1]}(x,y), T_{t}^{[2]}(x,y)$, $t>0$, $x,y\in X_+$ be related semigroups integral kernels,  respectively. Suppose that non-negative functions $h_1, h_2$ on $X_+$ are in $ RH_\8(\mu_+)$. Assume that,  for $j=1,2$, $t>0$, and $x,y\in X_+\setminus \Gamma$, we have:
\eq{\label{harm_two}
T_t^{[j]} h_j(x) = h_j(x),
}
\eq{\label{upper_gaussian}
0\leq T_{t}^{[j]}(x,y)\leq C \mu(B(x,\st))^{-1} e^{-\frac{d(x,y)^2}{c_1 t}},
}
\eq{\label{sem_N_est}
C^{-1} \mu_{h_j^2}(B(x,\st))^{-1} e^{-\frac{d(x,y)^2}{c_2 t}}\leq \frac{T_{t}^{[j]}(x,y)}{h_j(x)h_j(y)} \leq C \mu_{h_j^2}(B(x,\st))^{-1} e^{-\frac{d(x,y)^2}{c_3 t}}.
}
Notice that under the assumptions above Theorem~\ref{main_thm_1} describes $H^1_{L_j}(X_+)$ for $j=1,2$.

We say that a function $f$ on $X$ is even (or odd) if $f(\wt{x}) = f(x)$ (respectively, $f(\wt{x})=-f(x)$).  For each $f\in L^2(\mu)$ we decompose it into even and odd parts, $f=f_e+f_o$ getting $L^2(\mu) \simeq L^2(\mu_+)\oplus L^2(\mu_+)$, a direct sum of two Hilbert subspaces. We define $L$ as $L_1 \otimes L_2$, where $L_1$ acts on the even functions and $L_2$ acts on odd functions.  Such $L$ is selfadjoint and non-negative on $L^2(\mu)$. Moreover, $L$ is related to the semigroup $T_t$ whose integral kernel is given by:
\eq{\label{pol_skl}
T_t(x,y) =\frac{1}{2} 
\begin{cases}
    T^{[1]}_{t}(x,y) + T^{[2]}_{t}(x,y) & \text{for }x,y\in X_+\\
    T^{[1]}_{t}(x,\wt{y}) - T^{[2]}_{t}(x,\wt{y}) & \text{for }x, \wt{y}\in X_+\\
    T^{[1]}_{t}(\wt{x},y) - T^{[2]}_{t}(\wt{x},y) & \text{for }\wt{x},y\in X_+\\
    T^{[1]}_{t}(\wt{x},\wt{y}) + T^{[2]}_{t}(\wt{x},\wt{y}) & \text{for }\wt{x},\wt{y}\in X_+
\end{cases}.
}

We shall also assume that $T_t(x,y)$ satisfies \eqref{UG}. Notice that the functions
$$
\wh{h}_1(x) = \begin{cases}
    h_1(x) & \text{for }x\in X_+\\
    h_1(\wt{x}) & \text{for }x\in X_-
\end{cases},
\quad
\wh{h}_2(x) = \begin{cases}
    h_2(x) & \text{for }x\in X_+\\
    -h_2(\wt{x}) & \text{for }x\in X_-
\end{cases}
$$
are $L$-harmonic in the sense of \eqref{L-harm}.

\subsection{Hardy space related to $\wh{h}_1$ and $\wh{h}_2$.}

Our main goal of this section is to prove the following atomic characterization of $H^1_L(X)$.

\thm{thm_example_two}{
Let a space $X=X_+\cup X_-$ and a symmetry $\sigma$ satisfy the assumptions listed in Section \ref{ssec41}, including \eqref{metric}. Suppose that the conditions \eqref{harm_two}, \eqref{upper_gaussian}, \eqref{sem_N_est} hold for operators $L_1, L_2$ on $L^2(\mu_+)$ and functions $h_1, h_2\in RH_\infty(\mu_+)$, see \eqref{RHinfty}. Let $L = L_1\otimes L_2$ be the operator defined above and assume that the related semigroup kernel satisfies the upper Gaussian estimates \eqref{UG}. Then the Hardy space $H^1_L(X)$ coincides with $H^1_{at}(\mu; \wh{h}_1, \wh{h}_2; 2)$ and there exists $C$ such that
$$C^{-1} \norm{f}_{H^1_L(X)} \leq \norm{f}_{H^1_{at}(\mu; \wh{h}_1,\wh{h}_2;2)}\leq C\norm{f}_{H^1_L(X)}.$$
}

\begin{proof}
The inclusion $H^1_L(X) \subseteq H^1_{at}(\mu, \wh{h}_1, \wh{h}_2, 2)$ follows directly from Proposition \ref{prop_skr}, since \eqref{UG} holds for $L$ and both: $\wh{h}_1$ and $\wh{h}_2$ are $L$-harmonic functions.

Now we prove that  $H^1_{at}(\mu; \wh{h}_1, \wh{h}_2; 2)  \subseteq H^1_L(X)$. By the same continuity argument as in the proof of Theorem \ref{main_thm_1} (Step 2b on page \pageref{page7}) it is enough to prove that there exists $C>0$ such that
\eq{\label{atom-bdd}
\norm{M_La}_{L^1(X)} \leq C
}
for every $(\mu; \wh{h}_1, \wh{h}_2; 2)$-atom $a$.

Assume then that an atom $a$ is such that $\supp \, a \subseteq B=B(y,r) \subseteq X$. Denote $\wt{B} = \sigma(B)$ and
$$a_e(x) = \frac{1}{2} \eee{a(x) + a(\wt{x})}, \qquad a_o(x) = \frac{1}{2} \eee{a(x) - a(\wt{x})}.$$
Notice that $a_e$ and $a_o$ are supported in $B\cup \wt{B}$ and
$$\norm{M_L a}_{L^1(\mu)} \leq \norm{M_L a_e}_{L^1(\mu)}+\norm{M_L a_o}_{L^1(\mu)}.$$
First we shall estimate $\norm{M_L a_e}_{L^1(\mu)}$. Since $a_e, \wh{h}_1$ are even, $a_o$ is odd, and $\int a \wh{h}_1 \, d\mu(x) = 0$, we obtain
$$\int_{X_+} a_e h_1 \, d\mu = \frac{1}{2} \int_X a_e \wh{h}_1\, d\mu = \frac{1}{2} \int_X (a_e+a_o) \wh{h}_1\, d\mu=0. $$
Denote by $a_{e,+}$ the restriction of $a_e$ to $X_+$. Observe that, by \eqref{metric}, $a_{e,+}$ is supported on $B_+ :=B(z,r)\cap X_+ = (B\cup \wt{B})\cap X_+$, where $z=y$ if $y\in X_+$ or $z=\wt{y}$ otherwise. In addition, $\mu(B_+) \simeq \mu(B)$. Note that $T_t a_e(x)$ is an even function and $T_t a_e(x)=T_t^{[1]} a_{e,+}(x)$ for $x\in X_+$, so it is enough to show that
$$\norm{M_{L_1} a_{e,+}}_{L^1(\mu_+)} \leq C.$$
The last estimate is a consequence of Theorem \ref{main_thm_1}, since $a_{e,+}$ is an $(\mu;h_1;2)$-atom related to the ball $B_+\subseteq X_+$, see Definition~\ref{def_atoms}. 
The estimate for $\norm{M_L a_o}_{L^1(\mu)}$ is similar. We skip the details.

\end{proof}

Using Theorem \ref{thm_duality_general}, we  describe the dual space as follows.
\cor{cor_BMO_2}{
Assume that $h_1,h_2\in L^2_{loc}(\mu_+)$. Then the space $BMO(\mu;\wh{h}_1,\wh{h}_2;2)$ is dual to $H^1_{at}(\mu;\wh{h}_1,\wh{h}_2;2)$.
}

\section{Applications}
\label{sec5}
In this section we present exemplary applications of our theory, i.e. Theorems \ref{main_thm_1} and \ref{thm_example_two}.

\subsection{Operators with one harmonic function.}\label{sec5.0}

It is easy to check that all of the several examples discussed in \cite[Section  6]{Preisner_Sikora_Yan} satisfy the assumptions of Theorem \ref{main_thm_1}, including the condition $h\in RH_\8(\mu)$. Notice that in many cases the related harmonic function is unbounded either from above or from below. We do not list all examples here and, instead, we just quote one application.

{\it Exterior domain outside a bounded convex $C^{1,1}$ set.}
Assume that $\Omega\subset {\mathbb R}^n$ is such that $\Omega^c$
is convex, compact, and its boundary is locally a $C^{1,1}$ function.
We consider the Laplace operator $\Delta_\Omega$ on $\Omega$ with the Dirichlet boundary condition and the corresponding Hardy space $H^1_{\Delta_\Omega}(\Omega)$. 

\thm{thm_Dir2}{
Assume that for $n \ge 3$, $\Omega\subset {\mathbb R}^n$ is a domain that is an  exterior of a $C^{1,1}$ compact convex set.
 Then there exists a function $h:\Omega \to (0,\8)$, such that
$$h(x) \simeq \min(1, \mathrm{dist}(x, \Omega^c))$$
and the Hardy space $H^1_{\Delta_\Omega}(\Omega)$ coincides with $H^1_{at}(\mu; h; 2)$, where $\mu$ is the Lebesgue measure on~$\Omega$.
}
The estimates for the heat kernel were given in  \cite{Zhang_Bull_London_2003}, see also \cite[Section 6.1.2]{Preisner_Sikora_Yan}. It is not difficult to check that for $\Omega$ as in the theorem above we have $h\in RH_\8(\mu)$. Hence Theorem \ref{thm_Dir2} is a~direct consequence of Theorem \ref{main_thm_1}. Also, Theorem \ref{thm_duality} describes the dual space as $BMO(\mu;h;2)$. For other results concerning Hardy spaces on domains on $\RR^n$ we refers the reader to \cite{Auscher_Russ, Chang_Krantz_Stein} and references therein.

\subsection{Symmetric manifolds with ends}
\label{sec5.1}
In this section we consider a Riemmanian manifold with ends of the form $M = \Rn \# \Rn$. The results for this context were stated in Section \ref{sec13} and here we provide some details and precise assumptions. We shortly discuss the proof of Theorem~\ref{thm_manifold}. By definition, the manifold $\Rn \# \Rn$ has a central compact part $K$ such that $M \setminus K$ is a disjoint sum of two ends $E_+$ and $E_-$, where $E_+=E_-=\Rn \setminus B(0,r)$ and $r>0$. Let $\mu$ be a~measure on $M$ that has a smooth positive density and coincides with the Lebesgue's measure on $E_+$ and $E_-$. On the other hand we suppose that $M= M_+ \cup M_-$, where $E_+ \subseteq M_
+$ and $E_- \subseteq M_
-$. The halves $M_+$ and $M_-$ are Riemmanian manifolds  with boundaries and the set $\Gamma :=M_+\cap M_- = \mathrm{bd} \,M_+ = \mathrm{bd}\, M_- \subseteq K$ is a compact smooth sub-manifold of order $n-1$, i.e. $\mu(\Gamma)=0$.  Moreover, we assume that there exists a symmetry $\sigma : M \to M$ such that: $\sigma^2 = \mathrm{id}$, $\sigma$ preserves both the metric and the measure, $\sigma(x)=x$ for $x\in \Gamma$, and $\sigma(M_+)=M_-$. Notice that in this context the property \eqref{metric} holds.

Denote by $L_M$ the Laplace-Beltrami operator on $M$. It is known that the related heat kernel satisfies \eqref{UG}, c.f. \cite{Grigoryan_Saloff-Coste, Carron_Coulhon_Hassell}.

For a moment consider the half-space $M_+$ described above. Let $L_N$ and $L_D$ be the Laplace-Beltrami operators on $M_+$ with either Neumann or Dirichlet boundary conditions on $\Gamma$, respectively. We assume that $M_+$ is an inner uniform domain in the sense of \cite{Gyrya_Saloff-Coste}. By \cite[Theorem 3.10]{Gyrya_Saloff-Coste}, the heat kernel $T_{t,N}(x,y)$ corresponding to the semigroup $T_{t,N} = \exp(-tL_N)$ satisfy
	$$C^{-1} \mu(B(x,\st))^{-1} \exp\eee{-\frac{\rho(x,y)^2}{c_1 t}}\leq T_{t,N}(x,y) \leq C \mu(B(x,\st))^{-1} \exp\eee{-\frac{\rho(x,y)^2}{c_2 t}}$$
for $x,y \in M_+$ and $t>0$. Moreover, by \cite[Section 4 and Corollary 5.10]{Gyrya_Saloff-Coste} there exists a function $h_D:M_1 \to [0,\8)$ such that:
	$$C^{-1} \mu_{h_D^2}(B(x,\st))^{-1} \exp\eee{-\frac{\rho(x,y)^2}{c_1 t}}\leq \frac{T_{t,D}(x,y)}{h_D(x)h_D(y)} \leq C \mu_{h_D^2}(B(x,\st))^{-1} \exp\eee{-\frac{d(x,y)^2}{c_2 t}},$$
 where $t>0$ and $x,y\in M_+ \setminus \Gamma$. The function $h_D$ is continuous, positive on $M_+\setminus \Gamma$, and vanishes on $\Gamma$. Similarly as in Theorem \ref{thm_Dir2}, the function $h_D$ is in $RH_\8(\mu)$.

Observe that, thanks to symmetry, the analysis of $T_t$ on even functions on $M$ is equivalent to the analysis of $T_{t,N}$ on $M_+$, compare Section \ref{ssec41}. Similarly, the analysis of $T_t$ on odd functions on $M$ is equivalent to the analysis of $T_{t,D}$ on $M_+$. The kernels $T_t(x,y)$, $T_{t,N}(x,y)$ and $T_{t,D}(x,y)$ are related as in \eqref{pol_skl} with $T_t^{[1]} = T_{t,N}$ and $T_t^{[2]} = T_{t,D}$. On $M
_+$ the related harmonic functions are $h_1 \equiv 1$  and $h_2 = h_D$, respectively. Let $\wh{h}_1\equiv 1$ (the even extension of $h_1$ to $M$) and denote by $\wh{h}_2$ the odd extension of $h_D$ to $M$.

\begin{proof}
    [Proof of Theorem \ref{thm_manifold}]
    It is enough to observe that all the assumptions of Theorem \ref{thm_example_two} are satisfied with $T_t^{[1]} =T_{t,N}$ and $T_t^{[2]} = T_{t,D}$. Therefore $\wh{h}_1\equiv 1$ and $\wh{h}_2$ are $L$-harmonic. Finally observe that the functions $h_+$ and $h_-$ defined in Section \ref{sec13} are linear combinations of $\wh{h}_1$ and $\wh{h}_2$, i.e.
    $$h_+ = \frac{1}{2}\eee{\wh{h}_1+\wh{h}_2}, \qquad h_- = \frac{1}{2}\eee{\wh{h}_1-\wh{h}_2}.$$
\end{proof}

\rem{rem_rem}{
We intend to study the Hardy spaces on manifolds with ends without 
the symmetry assumption in another project.
}

{\bf Acknowledgements.}
We would like to thank Lixin Yan for reading our manuscript carefully and making valuable suggestions. 

{\bf Funding.}
M.P. and A.S. were  partly supported by Australian Research Council (ARC) Discovery Grants  DP200101065.


\begin{thebibliography}{10}

\bibitem{AMM}
P.~Auscher, A.~McIntosh, and A.~J. Morris, \emph{Calder\'{o}n reproducing
  formulas and applications to {H}ardy spaces}, Rev. Mat. Iberoam. \textbf{31}
  (2015), no.~3, 865--900.

\bibitem{Auscher_Russ}
P.~Auscher and E.~Russ, \emph{Hardy spaces and divergence operators on strongly
  {L}ipschitz domains of {$\Bbb R^n$}}, J.~Funct. Anal. \textbf{201} (2003),
  no.~1, 148--184.

\bibitem{Sikora_Bailey}
J.~Bailey and A.~Sikora, \emph{Vertical and horizontal square functions on a
  class of non-doubling manifolds}, J. Differential Equations \textbf{358}
  (2023), 41--102.

\bibitem{Carron_Coulhon_Hassell}
G.~Carron, T.~Coulhon, and A.~Hassell, \emph{Riesz transform and
  {$L^p$}-cohomology for manifolds with {E}uclidean ends}, Duke Math. J.
  \textbf{133} (2006), no.~1, 59--93.

\bibitem{Chang_Krantz_Stein}
D.-C. Chang, S.G. Krantz, and E.~M. Stein, \emph{{$H^p$} theory on a smooth
  domain in {${\bf R}^N$} and elliptic boundary value problems}, J. Funct.
  Anal. \textbf{114} (1993), no.~2, 286--347.

\bibitem{Coifman_Studia}
R.R. Coifman, \emph{A real variable characterization of {$H^{p}$}}, Studia
  Math. \textbf{51} (1974), 269--274.

\bibitem{Davies2}
E.~B. Davies, \emph{Non-{G}aussian aspects of heat kernel behaviour}, J. London
  Math. Soc. (2) \textbf{55} (1997), no.~1, 105--125.

\bibitem{Duong_Yan_2005}
X.~T. Duong and L.~Yan, \emph{Duality of {H}ardy and {BMO} spaces associated
  with operators with heat kernel bounds}, J. Amer. Math. Soc. \textbf{18}
  (2005), no.~4, 943--973.

\bibitem{DP_Annali}
J.~Dziuba\'{n}ski and M.~Preisner, \emph{Hardy spaces for semigroups with
  {G}aussian bounds}, Ann. Mat. Pura Appl. (4) \textbf{197} (2018), no.~3,
  965--987.

\bibitem{DZ_Revista2}
J.~Dziuba\'nski and J.~Zienkiewicz, \emph{On {H}ardy spaces associated with
  certain {S}chr\"odinger operators in dimension 2}, Rev. Mat. Iberoam.
  \textbf{28} (2012), no.~4, 1035--1060.

\bibitem{DZ_JFAA}
\bysame, \emph{On isomorphisms of {H}ardy spaces associated with
  {S}chr\"odinger operators}, J. Fourier Anal. Appl. \textbf{19} (2013), no.~3,
  447--456.

\bibitem{G_I_S-C2}
A.~Grigor'yan, S.~Ishiwata, and L.~Saloff-Coste, \emph{Geometric analysis on
  manifolds with ends}, Analysis and partial differential equations on
  manifolds, fractals and graphs, Adv. Anal. Geom., vol.~3, De Gruyter, Berlin,
  [2021] \copyright 2021, pp.~325--343.

\bibitem{G_I_S-C1}
\bysame, \emph{Poincar\'{e} constant on manifolds with ends}, Proc. Lond. Math.
  Soc. (3) \textbf{126} (2023), no.~6, 1961--2012.

\bibitem{Grigoryan_Saloff-Coste}
A.~Grigor'yan and L.~Saloff-Coste, \emph{Heat kernel on manifolds with ends},
  Ann. Inst. Fourier (Grenoble) \textbf{59} (2009), no.~5, 1917--1997.

\bibitem{Gyrya_Saloff-Coste}
P.~Gyrya and L.~Saloff-Coste, \emph{Neumann and {D}irichlet heat kernels in
  inner uniform domains}, Ast\'{e}risque (2011), no.~336, viii+144.

\bibitem{Sikora_Nix}
A.~Hassell, D.~Nix, and A.~Sikora, \emph{{Riesz transforms on a class of
  non-doubling manifolds {II}}}, arXiv e-prints (2019), arXiv:1912.06405, to
  appear in Indiana Univ. Math. J.

\bibitem{Sikora_Hassell}
A.~Hassell and A.~Sikora, \emph{Riesz transforms on a class of non-doubling
  manifolds}, Comm. Partial Differential Equations \textbf{44} (2019), no.~11,
  1072--1099.

\bibitem{Dangyang}
Dangyang He, \emph{{Endpoint Estimates For Riesz Transform And Hardy-Hilbert
  Type Inequalities}}, arXiv e-prints (2023), arXiv:2302.13739.

\bibitem{Hofmann_Memoirs}
S.~Hofmann, G.~Lu, D.~Mitrea, M.~Mitrea, and L.~Yan, \emph{Hardy spaces
  associated to non-negative self-adjoint operators satisfying
  {D}avies-{G}affney estimates}, Mem. Amer. Math. Soc. \textbf{214} (2011),
  no.~1007, vi+78.

\bibitem{Hofmann_Mayboroda_McIntosh}
S.~Hofmann, S.~Mayboroda, and A.~McIntosh, \emph{Second order elliptic
  operators with complex bounded measurable coefficients in {$L^p$}, {S}obolev
  and {H}ardy spaces}, Ann. Sci. \'{E}c. Norm. Sup\'{e}r. (4) \textbf{44}
  (2011), no.~5, 723--800.

\bibitem{Latter_Studia}
R.H. Latter, \emph{A characterization of {$H^{p}({\bf R}^{n})$} in terms of
  atoms}, Studia Math. \textbf{62} (1978), no.~1, 93--101.

\bibitem{Paluszynski_Zienkiewicz}
M.~Paluszynski and J.~Zienkiewicz, \emph{A remark on atomic decompositions of
  martingale {H}ardy's spaces}, J. Geom. Anal. \textbf{31} (2021), no.~9,
  8866--8878.

\bibitem{Preisner_Sikora_Yan}
M.~Preisner, A.~Sikora, and L.~Yan, \emph{Hardy spaces meet harmonic weights},
  Trans. Amer. Math. Soc. \textbf{375} (2022), no.~9, 6417--6451.

\bibitem{Song-Yan-2017}
L.~Song and L.~Yan, \emph{Maximal function characterizations for {H}ardy spaces
  associated with nonnegative self-adjoint operators on spaces of homogeneous
  type}, J. Evol. Equ. \textbf{18} (2018), no.~1, 221--243.

\bibitem{Song_Yan_2018}
\bysame, \emph{Maximal function characterizations for {H}ardy spaces associated
  with nonnegative self-adjoint operators on spaces of homogeneous type}, J.
  Evol. Equ. \textbf{18} (2018), no.~1, 221--243.

\bibitem{Stein}
E.~M. Stein, \emph{Harmonic analysis: real-variable methods, orthogonality, and
  oscillatory integrals}, Princeton Mathematical Series, vol.~43, Princeton
  University Press, Princeton, NJ, 1993, With the assistance of Timothy S.
  Murphy, Monographs in Harmonic Analysis, III.

\bibitem{Zhang_Bull_London_2003}
Qi~S. Zhang, \emph{A sharp comparison result concerning {S}chr\"odinger heat
  kernels}, Bull. London Math. Soc. \textbf{35} (2003), no.~4, 461--472.

\end{thebibliography}

\def\cprime{$'$}
\providecommand{\bysame}{\leavevmode\hbox to3em{\hrulefill}\thinspace}
\providecommand{\MR}{\relax\ifhmode\unskip\space\fi MR }
\providecommand{\MRhref}[2]{%
  \href{http://www.ams.org/mathscinet-getitem?mr=#1}{#2}
}
\providecommand{\href}[2]{#2}

\end{document}